\documentclass[manuscript]{acmart}
\usepackage[T1]{fontenc}
\usepackage[utf8x]{inputenc}
\usepackage{siunitx}
\usepackage{tikz}
\usepackage{natbib}
\usepackage{hyperref}
\usepackage{amsmath}
\usetikzlibrary{decorations.pathmorphing}
\usetikzlibrary{arrows.meta}

\graphicspath{{./gfx/}{./paper/gfx/}}

\numberwithin{equation}{section}

\AtBeginDocument{%
  }

\setcopyright{none}
\copyrightyear{2026}
\acmYear{2026}
\acmDOI{XXXXXXX.XXXXXXX}

\acmJournal{TOMS}
\acmVolume{0}
\acmNumber{0}
\acmArticle{0}
\acmMonth{0}

\title{Multi-Solver Coupling for Parallel Adaptive Multi-Physics Simulations with Trixi.jl and deal.II}
\author{Vivienne Ehlert}
\orcid{https://orcid.org/0000-0001-6586-2814}
\email{vivienne.ehlert@uni-a.de}
\affiliation{%
   \institution{University of Augsburg}
   \department{Institute of Mathematics, Centre for Advanced Analytics and Predictive Sciences}
   \city{86159 Augsburg}
   \country{Germany}
}
\author{Gregor Gassner}
\orcid{0000-0002-1752-1158}
\email{ggassner@uni-koeln.de}
\affiliation{%
   \institution{University of Cologne}
   \department{Department of Mathematics and Computer Science}
   \city{50931 Köln}
   \country{Germany}
}
\author{Martin Kronbichler}
\orcid{0000-0001-8406-835X}
\email{martin.kronbichler@rub.de}
\affiliation{%
   \institution{Ruhr University Bochum}
   \department{Department of Mathematics and Computer Science}
   \city{44801 Bochum}
   \country{Germany}
}
\author{Hendrik Ranocha}
\orcid{0000-0002-3456-2277}
\email{hendrik.ranocha@uni-mainz.de}
\affiliation{%
  \institution{Johannes Gutenberg University Mainz}
  \department{Institute of Mathematics}
  \city{55128 Mainz}
  \country{Germany}
}
\author{Michael Schlottke-Lakemper}
\orcid{0000-0002-3195-2536}
\email{michael.schlottke-lakemper@uni-a.de}
\affiliation{%
  \institution{University of Augsburg}
  \department{Institute of Mathematics, Centre for Advanced Analytics and Predictive Sciences}
  \city{86159 Augsburg}
  \country{Germany}
}
\acmJournal{TOMS}

\setcitestyle{authoryear} 

\newcommand{\fbf}{\mathbf{f}}
\newcommand{\kbf}{\mathbf{k}}
\newcommand{\sbf}{\mathbf{s}}
\newcommand{\ubf}{\mathbf{u}}
\newcommand{\vbf}{\mathbf{v}}
\newcommand{\xbf}{\mathbf{x}}

\newcommand{\nph}{\phantom{0}}

\newlength{\sedovw}
\newcommand{\sedovsplit}[1]{%
  \begin{tikzpicture}
    \node[anchor=south west,inner sep=0] (img) at (0,0)
      {\includegraphics[width=\sedovw]{#1}};
    \coordinate (splittop) at ([xshift=0.4058\sedovw]img.north west);
    \coordinate (splitbot) at ([xshift=0.4058\sedovw]img.south west);
    \draw[white,line width=1.2pt] (splittop) -- (splitbot);
    \node[anchor=north,inner sep=2pt,rounded corners=2pt,
          fill=white,fill opacity=0.85,text opacity=1]
      at ([xshift=0.2029\sedovw,yshift=-4pt]img.north west) {\footnotesize pressure};
    \node[anchor=north,inner sep=2pt,rounded corners=2pt,
          fill=white,fill opacity=0.85,text opacity=1]
      at ([xshift=0.6088\sedovw,yshift=-4pt]img.north west) {\footnotesize density};
  \end{tikzpicture}}

\begin{document}

\begin{abstract}
   Many standalone frameworks for numerical solvers have been developed to tackle the simulation of specific single- or multi-physics problems. For solving coupled problems, common approaches are to extend existing solvers, to develop an entirely new solver or to couple two existing solvers by using the interface provided by a coupling library or framework.
   However, to the best of our knowledge, there is not yet a framework that allows the easy and direct development of coupled solvers.
   In this work, we prototype a portable reproducible cross-language framework for the development of coupled parallel adaptive solvers using Trixi.jl and deal.II for the numerical simulation of coupled multi-physics problems. Currently, this is tightly entangled with an example coupled solver. We show its usability by developing a partitioned strongly-coupled multi-physics solver for the dynamics of Newtonian self-gravitational gases. For the coupled solver, we validate the expected order of convergence, physical sensibility of the results and mesh adaptivity. Finally, we investigate its parallel scaling. A publicly accessible reproducibility repository for the numerical results and code is available.
\end{abstract}




\maketitle

\section{Introduction} \label{sc:intro}

Many phenomena occurring in physics or engineering can best be described by combining multiple single-physics models into a jointly coupled multi-physics system \cite{holz2025}. Examples for such multi-physics systems are fluid flow against elastic materials, often called fluid-structure interaction, the flow of a fluid with particles suspended in it, or the interaction of electrically conducting fluid with magnetic fields. These multi-physics systems can be distinguished by the kind of coupling between the different constituting physical systems. Fluid-structure interaction is an example of systems in which the single-physics systems are coupled along abstract lower-dimensional interfaces along which they affect each other, while the interaction of electrically conducting fluids with magnetic fields is an example of systems in which the single-physics systems occupy the same physical space and interact with each other at every point of the computational domain. These interactions can be, for example, the exchange of energy and momentum between the coupled single-physics systems \cite{keyes2013}.

Mathematically, the multi-physics system can be modeled by sets of governing partial differential equations that may be of different types (e.g.~hyperbolic-elliptic, hyperbolic-parabolic or elliptic-parabolic), where each of the different types models different physical phenomena and crucially influences the way the coupled problem is solved \cite{keyes2013}. Hyperbolic partial differential equations can be advanced locally, as information travels according to a finite propagation speed, whereas in elliptic or parabolic partial differential equations the same information affects the whole computational domain immediately. In consequence, different solution strategies are required for different fields. Essentially, the way the set of governing equations is handled in an implementation in a solver can be distinguished either as being a unified, monolithic model or as being a partitioned, modular model. Both of these choices offer different possibilities in their implementations. Also implementations using each of these models can be either strongly or loosely coupled where strong couplings are simpler to achieve when using a unified model \cite{holz2025}.

For some split single-physics equations it is necessary to use specialized or structure-preserving methods to get stable, consistent and convergent solvers. This can be easier to implement consistently following a partitioned strategy, whereas a unified strategy offers an easier implementation of strongly coupled solvers. These kind of specialized solvers include entropy-stable time integration methods, implicit solvers, exterior calculus compliant spatial discretizations along with appropriate time integration schemes as well as Lie group methods.

In this work, we exemplify the proposed abstract coupling strategy by the dynamics of non-relativistic self-gravitating gases, which can be modeled by the compressible Euler equation for hydrodynamics and a Poisson equation for the self-gravitational potential, both coupled in bulk via linear source terms depending on each others state variables. For this problem, many solution strategies have been published previously using a variety of different approaches. In many works \cite{almgren2010,tomida2023,fukushima2026,flash}, multigrid methods for the Poisson operator are used to determine the self-gravitational potential. \citet{flash} propose further algorithmic options like a multi-pole expansion for near-spherical mass distributions. Some works also lift the Poisson equation to a system of hyperbolic equations which is then pseudo-time stepped until a tolerance is achieved to regain the solution to the Poisson equation, e.g.~in \cite{hirai2016,schlottkelakemper2021}. In \cite{hubber2018} a Lagrangian approach is used in which smoothed-particle hydrodynamics are used. The efficiency and hence the choice of the solver for the Newtonian gravitational potential is dependent on the distribution of density in the computational domain. Multigrid methods work well for spread-out and general mass densities, while for single large nearly spherical mass densities multi-pole expansions using spherical harmonics can be more efficient \cite{couch2013}. For a general performance assessment of multigrid and multipole methods, see, e.g.,~\cite{Gholami2016}.

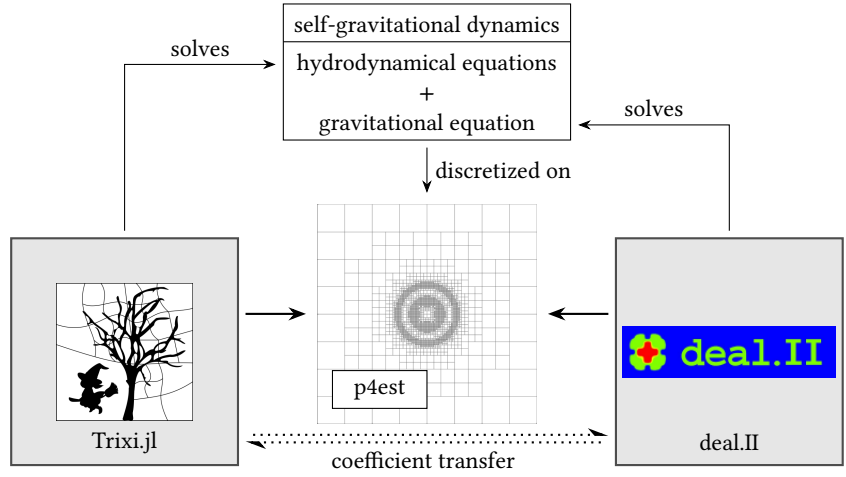
\begin{figure}
    \centering
    \begin{tikzpicture}
      \filldraw[thick,color=white!30!black,fill=white!90!black] (0,0)--(3,0)--(3,3)--(0,3)--cycle;
      \node at (1.5,0.3) {Trixi.jl};
      \node at (1.5,1.5) {\includegraphics[scale=0.12]{{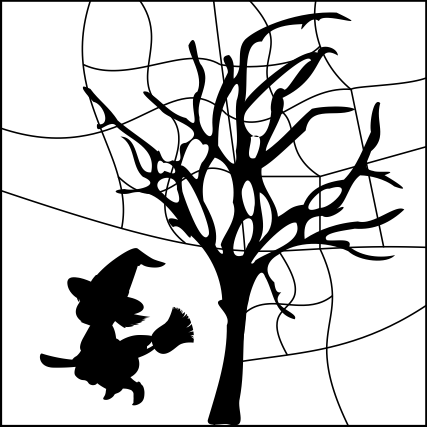}}};      
      \filldraw[thick,color=white!30!black,fill=white!90!black] (8,0)--(11,0)--(11,3)--(8,3)--cycle;
      \node at (9.5,0.3) {deal.II};
      \node at (9.5,1.5) {\includegraphics[scale=0.4]{{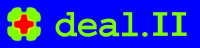}}};
      \node at (5.5,2.0) {\includegraphics[scale=0.8]{{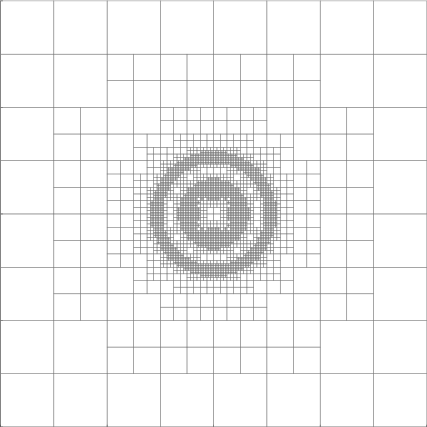}}};
      \filldraw[fill=white] (4.25,0.75) -- (5.50,0.75) -- (5.50,1.25) -- (4.25,1.25) -- cycle;
      \node at (4.875,1.0) {p4est};
      \draw[thick,-Stealth] (3.1,2.0) -- (3.9,2.0);
      \draw[thick,Stealth-] (7.1,2.0) -- (7.9,2.0);
      \draw[thick,dotted,arrows={-Stealth[harpoon]}] (3.2,0.4) -- (7.9,0.4);
      \draw[thick,dotted,arrows={-Stealth[harpoon]}] (7.8,0.3) -- (3.1,0.3) node[midway,below] {coefficient transfer};

      \draw (3.6,6.1) -- (3.6,4.3) -- (7.4,4.3) -- (7.4,6.1) -- cycle;
      \draw (3.6,5.6) -- (7.4,5.6);
      \node at (5.5,5.8) {self-gravitational dynamics};
      \node at (5.5,5.3) {hydrodynamical equations};
      \node at (5.5,4.9) {{\LARGE + }};
      \node at (5.5,4.5) {gravitational equation};

      \draw[-Stealth] (5.5,4.2) -- (5.5,3.6) node[right,midway] {discretized on};
      \draw[-Stealth] (1.5,3.1) -- (1.5,5.3) -- (3.5,5.3) node [above,sloped,align=center,midway] {solves};
      \draw[-Stealth] (9.5,3.1) -- (9.5,4.5) -- (7.5,4.5) node [above,sloped,align=center,midway] {solves};

    \end{tikzpicture}
    \Description{This figure illustrates the abstract coupling idea for the coupled solver we develop in this work. As we consider here the dynamics of self-gravitational gases of which the governing equations consist of the hydrodynamical equations and a gravitational equation coupled via source terms, of which we want to solve the hydrodynamical equations using Trixi.jl, a framework for numerical simulations of conservation laws, and the gravitational equation we want to solve using deal.II, a finite element library, there are arrows pointing from Trixi.jl and deal.II to the respective equations. The governing equations are discretized on a common mesh of which the spatial discretization is handled by p4est which is shared by Trixi.jl and deal.II. During the simulation, however, the partitioned solvers from Trixi.jl and deal.II have to communicate the values of their current coefficients with one another.}
    \caption{An illustration showing the abstract coupling procedure we follow in this work.}
    \label{fig:intro:coupling}
\end{figure}

In this work, we develop a partitioned solver for the dynamics of non-relativistic self-gravitational gases, which essentially follows the approach taken in \cite{schlottkelakemper2021}. To the hydrodynamical equations we apply the semi-discrete operator of the discontinuous Galerkin spectral element method and to the Newtonian self-gravitational equation we apply a multigrid solver. In order to reflect the general coupling setup with specialized software strategies for each of the two fields, we use two different numerical solver frameworks, Trixi.jl and deal.II. As a result, separate solvers are applied to the split equations while coupling them strongly in the time integration. The strong coupling is achieved by solving for the self-gravitational potential before each of the Runge--Kutta stages during time evolution of the hydrodynamical equations. The governing equations are spatially discretized on a common mesh handled by {\ttfamily p4est} \cite{BursteddeWilcoxGhattas2011}. Fig.~\ref{fig:intro:coupling} illustrates the chosen approach. 
The developed coupled solver is a prototype towards a modular framework for developing arbitrary coupled solvers using the Trixi.jl and deal.II infrastructure. The solver also supports parallel and \(h\)-adaptive execution. At the same time, the setup is still easy enough to understand for researchers to experiment with the ingredients and extend them further.


This work is organized in the following way. We start with brief introductions of the considered multi-physics problem at hand and its governing equations. Then the methods to be used for the numerical solution as well as the considered software frameworks in which these methods are implemented are presented. Next, we describe the bulk coupling implementation via a partitioned approach and the software steps that ensure reproducibility and portability for the coupled solver. Afterwards, we perform numerical experiments to validate and benchmark our solver. First, we consider a manufactured convergence test and the Jeans instability as tests to validate whether our method achieves the expected convergence rate and to validate it against analytical energy profiles available for the Jeans instability. We conclude this work after using our coupled solver to simulate cylindrical and spherical Sedov blast waves and benchmark our solver on simulations of spherical Sedov blast waves.

\section{Governing equations and numerical methods} \label{sc:govnum}

In this section, we first introduce the multi-physics problem that is considered in this work as a test problem. Afterwards we outline the discontinuous Galerkin spectral element method (DGSEM) used for the spatial discretization of the governing equations. Then we will give a short remark about adaptive mesh refinement and \(h\)-adaptivity of finite element discretizations. This is followed by a short description of multigrid methods followed by a brief introduction of the numerical simulation frameworks with which we develop a coupled solver in this work.

\subsection{Problem setting} \label{sc:govnum:setting}

We consider the dynamics of non-relativistic self-gravitating gases in two and three spatial dimensions as test multi-physics problem for which we want to develop a coupled solver combining Trixi.jl and deal.II. The state variables are the density of the gas \(\rho\), its momentum \(\rho\,v\), its total energy \(E\) and the self-gravitational potential \(\phi\). The governing equations can be split into the hyperbolic compressible Euler equations with source terms given by the self-gravitational potential in their momenta and energy and the elliptic self-gravitational equation with the density of the gas as a source \cite{chandrasekhar1961}. The compressible Euler equations in arbitrary dimensions can be expressed in conserved quantities as
\begin{equation}
    \frac\partial{\partial t}\begin{bmatrix} \rho \\[2pt] \rho\,v \\[2pt] E\end{bmatrix} + \nabla\cdot\begin{bmatrix} \rho\,v \\[2pt] \rho\,v \otimes v + p\mathbf{I} \\[2pt] (E+p)\,v\end{bmatrix} = \begin{bmatrix} 0 \\[2pt] -\rho\,\nabla\phi \\[2pt] -(v\cdot\nabla\phi)\,\rho\end{bmatrix}. \label{eq:govnum:compeuler}
\end{equation}
The pressure \(p\) can be determined from the state variables using the ideal gas law, with the heat capacity ratio \(\gamma\) for the gas, by
\begin{equation}
    p = (\gamma - 1)\left(E -\frac\rho2\|v\|^2\right). \label{eq:govnum:pressure}
\end{equation}
The equation for the non-relativistic self-gravitational potential is modeled by a Poisson equation with right-hand side given by the density of the gas:
\begin{equation}
    \Delta\,\phi = 4\pi\,G\,\rho, \label{eq:govnum:selfgrav}
\end{equation}
where \(G\) is the gravitational constant.

The source terms in Eqs.~\eqref{eq:govnum:compeuler} and \eqref{eq:govnum:selfgrav} give rise to a coupling in volume, as at each point of the computational domain they influence each other by the value of their state variables.

\subsection{Discontinuous Galerkin Spectral-Element Method} \label{sc:govnum:dgsem}

The discontinuous Galerkin spectral-element method (DGSEM) is a method used to discretize partial differential equations in space, see, e.g., \cite{kopriva2009} for a comprehensive description.
We assume a hyperbolic conservation law defined on a domain \(\Omega\) given by
\begin{equation}
    u_t(x) + \nabla \cdot f(u(x)) = s(u(x)),\qquad\forall x\in\Omega.
\end{equation}
Here \(u\in\mathbb{R}^m\), where \(m\) is the number of equations, is the state vector of the conserved variables in the dynamical system. The vector \(f = (f_1, \ldots, f_d)^\top\) denotes the physical fluxes and \(s\) denotes the source terms, which might also vanish. \(\Omega\) is a \(d\)-dimensional manifold which gets subdivided into \(k\) non-overlapping elements. Each element is reparameterized to the reference element \(E=[-1,1]^d\). The subdivided and reparameterized partial differential equation on the reference element, in reference coordinates \(\xi\), is then
\begin{equation}
    J\,u_t + \nabla_\xi \cdot f = J\,s(u),\label{eq:govnum:dgsem:reparametrized}
\end{equation}
where \(J\) is the determinant of the Jacobian of the parameterization, which in case of Cartesian elements on a hypercube has the value \(J=(h/2)^d\), where \(h\) is the size of the element.

The starting point of the DGSEM is the weak form of hyperbolic conservation laws. Multiplying Eq.~\eqref{eq:govnum:dgsem:reparametrized} with an appropriate test function, integrating over the reference element \(E\) and performing integration by parts, the weak form on each element can be expressed as
\begin{align}
    \int_E J\,u_t\,\varphi\;d\xi + \int_{\partial E} (f\cdot n)\,\varphi\;dS - \int_E f\cdot\nabla_\xi\varphi\,d\xi = \int_E J\,s\,\varphi\;d\xi\;,\label{eq:govnum:dgsem:weakform}
\end{align}
where \(n\) is the outer unit normal of the boundary. We expand the test functions in the space of polynomials of degree \(N\) as linear combinations of the tensor product of the Lagrange basis with respect to the \(N+1\) Legendre--Gauss--Lobatto (LGL) nodes \(\{x_k\}_{k=0}^N\) on \([-1,1]\):
\begin{align}
    \varphi(\xi) = \prod_{j=1}^d \sum_{k_j=0}^N \varphi_{k_j}\,\ell_{k_j}(\xi).
\end{align}
We insert this expansion into the weak form in Eq.~\eqref{eq:govnum:dgsem:weakform} for the test functions \(\varphi\) and analogously for the solution \(u\) and fluxes \(f\). By replacing the boundary fluxes \(f\cdot n\) by numerical fluxes \(f^\star_n\) between adjacent elements, i.e.~by numerical solutions to the Riemann problem on the boundaries of adjacent elements, one gets the weak form in the standard DGSEM formulation:
\begin{align}
    \int_{E} J\,u_t\,\varphi_i\;d\xi + \int_{\partial E} f^\star_n\,\varphi_i\;dS - \int_{E} f\cdot\nabla_\xi\varphi_i\,d\xi = \int_{E} J\,s\,\varphi_i\;d\xi\;,\quad\forall i\in \{0,\ldots,N\}^d\label{eq:govnum:dgsem:dgsem}
\end{align}
where the integrals are evaluated using LGL-quadrature. This leads to a collocation of interpolation and quadrature nodes. When actually inserting the expansions for flux and solution and requiring relation~\eqref{eq:govnum:dgsem:dgsem} to hold for each test function, the relation turns into a semi-discretization of the partial differential equation. Subsequently, a temporal discretization using an explicit Runge--Kutta-method or any other time integration method is made. For stable explicit time integration, the time step is chosen as
\begin{align}
    \Delta t = \frac{{\text{CFL}}}{N+1}\frac{h}{\Lambda},
\end{align}
where \(\Lambda=\max(|v| + c)\) with the speed of sound \(c=\sqrt{\gamma\,p/\rho}\) for the compressible Euler system.

Alternatively, the summation-by-parts (SBP) property of nodal discontinuous Galerkin schemes on LGL nodes can be used to obtain a split-form DGSEM approximation \cite{gassner2013}. The strong form of the conservation law in Eq.~\eqref{eq:govnum:dgsem:reparametrized} can be obtained by applying integration-by-parts again to Eq.~\eqref{eq:govnum:dgsem:dgsem}, an operation formally corresponding to a discrete summation-by-parts procedure when applying the LGL quadrature. A numerical volume flux for the flux derivative can be introduced, see \cite{carpenter2014,fisher2013}, to yield a split form discontinuous Galerkin approximation \cite{gassner2016,gassner2018}, allowing the use of entropy conserving or kinetic energy preserving symmetric two-point flux functions \cite{tadmor2003,ismail2009,chandrashekar2013,ranocha2018}. The numerical flux of \citeauthor{chandrashekar2013}, \cite{chandrashekar2013}, is used here.

For simulating the compressible Euler equations with strong discontinuities, a high-order shock capturing approach according to~\cite{hennemann2021} is utilized in the discontinuous Galerkin solver. In their approach, each element of the DG discretization is divided into \((N+1)^d\) subcells, on which a first-order finite volume method is used for obtaining a semi-discrete operator. Blending the finite volume based propagation with the tentative DG solution based on the energy content in the highest modes in each element yields an operator retaining the discrete entropic property.

For consistent treatment of local adaptive mesh refinement with hanging nodes, we use the mortar method \cite{Bernardi1994,kopriva1996,kopriva2002}. There, mortar surfaces are inserted at interfaces of coarse and refined elements. At non-conforming element interfaces the solution values are interpolated to the mortar, where the surfaces flux is calculated at conforming node locations. The resulting flux values are then discretely projected back to the non-conforming element interfaces. A similar approach is used for \(h\)-adaptivity of DG elements. For the refining of a cell, the solution on the original coarse element is interpolated to the LGL nodes on the refined elements, while for coarsening of cells the solution on the \(2^d\) refined elements is projected onto the coarse element. In \cite{bui-thanh2012} it is shown that the mortar method and the mesh refinement technique, which corresponds to a mortar element method in a higher-dimension, used here are fully conservative. In \cite{schlottke-lakemper2019}, more details for deriving the DGSEM on non-conforming hierarchical Cartesian meshes can be found. However, it should be noted that this mortar method introduces a two to one balancing along the interfaces. The description of this method is also easily lifted to higher dimensions.

\subsection{Multigrid for elliptic problems} \label{sc:govnum:mg}

The elliptic Poisson equation \eqref{eq:govnum:selfgrav} for the self-gravitational potential is solved with a multigrid approach.
Multigrid methods are efficient solution techniques for solving the linear systems arising from elliptic partial differential equations, such as the Poisson equation for the self-gravitational potential in Eq.~\eqref{eq:govnum:selfgrav}. The method combines simple iterative schemes, which are effective in reducing the high-frequency portion of the error, on a hierarchy of different grid resolutions \cite{trottenberg2001}. We employ three key ingredients to use as efficient a realization as possible \cite{kronbichler2018,Clevenger2021,munch2023}. Firstly, the hierarchy is expressed as a high-order DG scheme with the same nodes as the LGL basis in the hyperbolic solver component, where auxiliary-space arguments allow to embed the DG scheme into a continuous Galerkin representation of the same order~\cite{Antonietti2017}. Secondly, a combination of polynomial ($p$) and geometric ($h$) multigrid using the level transfer described in~\cite{munch2023} is used to construct the coarser problem sizes. On each level, a smoother based on a Chebyshev iteration of the point-Jacobi scheme is used, which strikes a good balance of cost and effectiveness for degrees $1\leq k \leq 8$ on Cartesian meshes~\cite{Adams2003,kronbichler2018}. As a last ingredient, we use matrix-free operator evaluation to implement fast matrix-vector products within the Chebyshev iteration, for residuals and the level transfer. Matrix-free evaluation does not store a global assembled matrix, but evaluates the finite-element integrals on the fly. This is beneficial for modern hardware, where computations are cheap compared to memory access, thus allowing to trade memory transfer for some additional computations. For the given quadrilateral and hexahedral meshes with tensor-product basis functions, sum factorization is utilized. A consistent Gauss--Legendre quadrature formula is used rather than the LGL quadrature in the DG field, because the cost for operator evaluation is often limited by the memory access also for matrix-free methods~\cite{Kronbichler2019}. The multigrid hierarchy is used as a preconditioner within a conjugate gradient algorithm to increase the robustness of the solver, leading to slightly lower iteration counts. Furthermore, the multigrid hierarchy is run in single-precision, with the outer CG solver done in double precision, to increase the solver throughput~\cite{Kronbichler2019}.

\subsection{Adaptive mesh refinement} \label{sc:govnum:amr}

Using a local truncation error estimator along with a threshold allows one to determine whether an element of the mesh should be refined into smaller cells or whether it should be coarsened \cite{berger1989}. However, adapting the mesh has to take the two-to-one balancing, which is due to the used mortar method, into account. Cells can only be coarsened if all their sibling cells, i.e.~the child nodes of their parent, are leaves, if the coarsen indicator on each cell is set and if they do not have to stay refined for balancing reasons. The usual algorithm for adaptation consists of first refining the marked cells along with ensuring that enough other cells are refined to satisfy the imposed balancing. Afterwards cells that can be coarsened without breaking the balancing can be coarsened. If the order of refining and coarsening were switched, one would potentially coarsen cells in the first step that need to be refined in the refinement step again to satisfy the balancing of the mesh.

The shock capturing indicator by \citeauthor{hennemann2021} \cite{hennemann2021} can also be used as a mesh refinement indicator by refining cells in which shocks occur and coarsening cells in which no shocks are determined. The two to one balancing imposed by the used mortar method yields then a smoothly refined mesh.

\subsection{Numerical solvers} \label{sc:govnum:numsolvers}

In this work, we couple Trixi.jl and deal.II for performing numerical simulations of the multi-physics system given by Eqs.~\eqref{eq:govnum:compeuler} and \eqref{eq:govnum:selfgrav}. We use Trixi.jl to solve Eq.~\eqref{eq:govnum:compeuler} and deal.II to solve Eq.~\eqref{eq:govnum:selfgrav}.

Trixi.jl \cite{schlottkelakemper2025} is an open source framework written in the Julia programming language~\cite{bezanson2017} for adaptive high-order numerical simulations of hyperbolic conservation laws implementing the discontinuous Galerkin spectral element method (DGSEM) described in subsection~\ref{sc:govnum:dgsem}. The development of this framework is based on a modular architecture, in which all of the components are only loosely coupled with one another. Hence, making it easy to extend or replace existing functionality.

The open source finite element library deal.II is written in the C++ programming language and contains abstractions for common algorithms of adaptive finite element methods for a broad variety of partial differential equations \cite{arndt2021,arndt2025}. In this work, we make use of the matrix-free and multigrid modules of the library as described in subsection~\ref{sc:govnum:mg}, providing an efficient solver for the non-relativistic self-gravitational potential.

\section{Coupling method and algorithms} \label{sc:algs}

In this section, we describe a coupling strategy of two numerical solvers, starting from an abstract description and the main steps to achieve the coupling operations. The concept is then exemplified by the realization and technical details of a coupling between the Julia-based solver Trixi.jl and the C++-based framework deal.II. We close this section with a brief description of how we made our coupled solver accessible, reproducible and portable across a selection of different platforms. The implementations in C++ and Julia of the algorithms and methods described in this section can be found in the following reproducibility repository hosted on GitHub: \cite{ehlert2026multisolverRepro}.

\subsection{Coupling} \label{sc:algs:coupling}

The equations of state of the dynamical system given in subsection~\ref{sc:govnum:setting} form a bulk-coupled problem. Depending on the solver approach, this can result in an increase of complexity in the implementation, as more information may need to be communicated between the fields compared to interface couplings over abstract, lower dimensional interfaces where different physical systems interact~\cite{preCICEv2}. This is especially true for partitioned approaches, where one considers the constituting equations separately compared to monolithic approaches, where the whole unified set of governing equations is considered at each point in the computational domain.

This work considers a partitioned approach, in order to use specialized solvers for each of the constituting equations. Hence, as the governing equations are the hyperbolic compressible Euler equations, Eq.~\eqref{eq:govnum:compeuler}, and the elliptic non-relativistic self-gravitational Poisson equation, Eq.~\eqref{eq:govnum:selfgrav}, we will solve Eq.~\eqref{eq:govnum:compeuler} with an efficient implementation of the DGSEM in Trixi.jl and Eq.~\eqref{eq:govnum:selfgrav} with an efficient multigrid method implemented using deal.II.

As a data exchange mechanism, our approach is to let both solvers agree on a common representation of bulk solution vectors using the same element description. The DGSEM approach used on the Trixi.jl side employs polynomials collocated in the Legendre--Gauss--Lobatto points, which can be matched by the elliptic solver with the element \texttt{FE\_DGQ<d>(degree)} of the deal.II library~\cite{arndt2021}. During the numerical simulation, the solvers communicate their state variables in this format. As the multigrid solver eventually employs continuous finite elements, it implements a mapping in the form of a projection of the coefficients of the discontinuous elements onto the continuous elements in the right-hand side of Eq.~\eqref{eq:govnum:selfgrav}, using a test function $\varphi$, in the usual finite element way,
\[
\int_E -4 \pi G \rho \varphi J \, d\xi.
\]
Since any $L^2$-conforming function $\rho$ can be paired with an $H^1$-conforming test function $\varphi$, the operation does not need additional measures to deal with the non-conformity. To realize the regularity gain of the Poisson solver, the continuous finite element space uses polynomials of degree $k+1$, if the DG solver employs polynomials of degree $k$. After the solution of the Poisson problem, the solution is projected back to the discontinuous function space. Due to the truncation from the space $\mathbb{Q}_{k+1,0}$ to $\mathbb{Q}_{k,-1}$ a small projection error arises.

A necessary ingredient for the data exchange to happen solely via solution vectors is that the two solvers need to share the same mesh to discretize the spatial computational domain. To enable a low-invasive coupling, each solver constructs its own mesh, starting from a mesh initially created on the Julia side of the implementation in Trixi.jl, which is then synchronized to the C++ side of the implementation in deal.II by means of the relevant refinement indicators. This is a small overhead, but justified for higher-order methods where the mesh storage itself only constitutes a small fraction of the overall memory consumption and computational effort. Compliance with the native data structures of each package eases the implementation. This strategy forces the solver to communicate all the changes made to the mesh in Trixi.jl in the process of adaptive mesh refinement. Note that the communication would need to happen in either case, even if there were a unique instance of the mesh being shared, as the other solver needs to be aware of changes in the interpretation of the coefficient vectors and incidence relations between elements of the mesh.

The mesh handling is performed by {\ttfamily p4est}, \cite{BursteddeWilcoxGhattas2011}, which is supported by both solvers, see also \cite{Bangerth2011}. This further simplifies interoperation as we can have a one-to-one loop through cells without additional geometric searches. Also, as {\ttfamily p4est} supports MPI parallel usage, and both Trixi.jl and deal.II support MPI parallel usage, our coupled solver can also be used in parallel with MPI. As the mesh distribution over MPI ranks is deterministic with this setup, the partitioned meshes of the solvers on the respective MPI ranks are guaranteed to coincide with each other.

\begin{figure}
    \centering
    \begin{tikzpicture}
        \draw[->] (-0.25,0) node[left] {time discretization} -- (8.5,0);
        \draw[->] (-0.25,1.5) node[left] {coupled solver} -- (8.5,1.5);
        \draw[->] (-0.25,3.0) node[left] {callbacks} -- (8.5,3.0);
        \filldraw[color=black,fill=white] (0,-0.25) -- (2,-0.25) -- (2,0.25) -- (0,.25) -- cycle;
        \filldraw[color=black,fill=white!90!black] (2,1.25) -- (3,1.25) -- (3,1.75) -- (2,1.75) -- cycle;
        \filldraw[color=black,fill=white] (3,-0.25) -- (4,-0.25) -- (4,0.25) -- (3,.25) -- cycle;
        \filldraw[color=black,fill=white] (4,2.75) -- (5,2.75) -- (5,3.25) --node[midway,below] {stage} (4,3.25) -- cycle;
        \filldraw[color=black,fill=white] (5,-0.25) -- (6,-0.25) -- (6,0.25) -- (5,.25) -- cycle;
        \filldraw[color=black,fill=white] (6,2.75) -- (7,2.75) -- (7,3.25) -- node[midway,below] {step} (6,3.25) -- cycle;
        \filldraw[color=black,fill=white] (7,-0.25) -- (8,-0.25) -- (8,0.25) -- (7,.25) -- cycle;
        \draw[-Stealth] (0,-.75) -- (0,-.25);
        \draw[-Stealth] (2,0.25) -- node[midway,sloped,above]{rhs!} (2,1.25);
        \draw[-Stealth] (3,1.25) -- (3,0.25);
        \draw[-Stealth] (4,0.25) -- node[midway,sloped,above]{callback!}(4,2.75);
        \draw[-Stealth] (5,2.75) -- (5,0.25);
        \draw[-Stealth] (6,0.25) -- node[midway,sloped,above]{callback!}(6,2.75);
        \draw[-Stealth] (7,2.75) -- (7,0.25);
        \draw[-Stealth] (8,-.25) -- (8,-.75);
        \draw[thick,densely dotted] (1.5,-0.5) -- (5.5,-0.5) -- (5.5,4.0) -- (1.5,4.0) -- cycle;
        \draw[thick] (1.5,3.6) -- (2.5,3.6) -- (2.6,3.7) -- (2.6,4.0) -- (1.5,4.0) -- cycle;
        \node at (2.0,3.8) {loop}; \node at (4.0,3.8) {per RK stage};
        \draw[thick,densely dotted] (0.5,-.75) -- (7.5,-.75) -- (7.5,4.75) -- (0.5,4.75) -- cycle;
        \draw[thick] (0.5,4.35) -- (1.5,4.35) -- (1.6,4.45) -- (1.6,4.75) -- (0.5,4.75) -- cycle;
        \node at (1.0,4.55) {loop}; \node at (4.0,4.55) {per time step};

    \end{tikzpicture}
    \Description{A pseudo-UML sequence diagram showing the time integration loop dissembled into lifelines for the time discretization, semi-discretization and callbacks.}
    \caption{A timeline showing the calling relations during time integration per time step, following the structure of general Trixi.jl solvers. The calling sequence in the coupled solver, whose activation box is highlighted here in grey, is detailed in Fig.~\ref{fig:algs:seq-rhs}.}
    \label{fig:algs:timeline}
\end{figure}
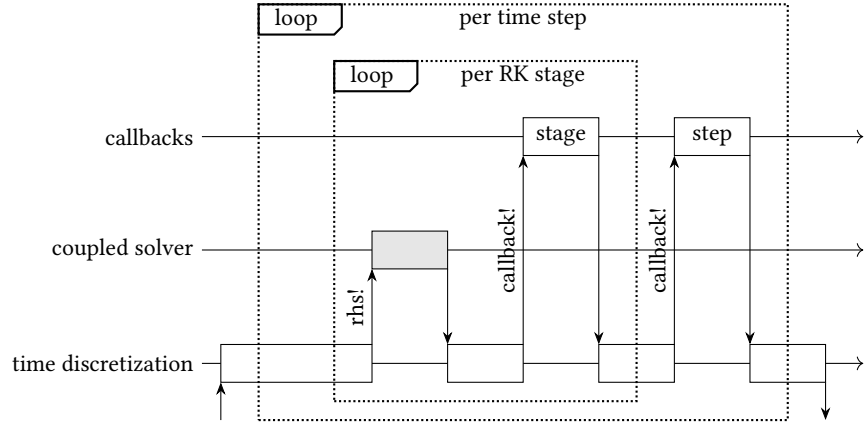

As only one of the governing equations~\eqref{eq:govnum:compeuler} and~\eqref{eq:govnum:selfgrav} contains a time derivative, the main control lies in the time integration loop for Eq.~\eqref{eq:govnum:compeuler}. It has been shown that arbitrary explicit Runge--Kutta schemes are appropriate for time integration of this coupled system \cite{schlottkelakemper2021}. However, for high-order numerical accuracy of the time integration scheme the self-gravitational potential has to be solved for in each of the Runge--Kutta stages. An illustration of the main time integration loop in Trixi.jl for this dynamical system is given in Fig.~\ref{fig:algs:timeline}. Do note, that the size of the activation boxes in the sequence diagrams do not correspond to time or work spent in the called subroutines. The figure also shows that the Runge--Kutta time step evaluation can be followed by several system callbacks for performing I/O, time step control or adaptive mesh refinement.

\subsection{Technical details} \label{sc:algs:technicals}


Since only the compressible Euler equations in Eq.~\eqref{eq:govnum:compeuler} carry a time derivative, the main control loop is the time integration loop for the Euler system. As the time integration loop in Trixi.jl is handled via semi-discretizations, the solver is implemented as a semi-discretization for the coupled equations containing two further semi-discretizations respectively for the hydrodynamical equations and the self-gravitational equation.

The time integration loop in Trixi.jl is implemented to be compatible with the OrdinaryDiffEq.jl package, which provides many different time integration schemes. To this end, the semi-discretization overloads a function to implement a right-hand-side evaluation according to the governing hydrodynamic equations~\eqref{eq:govnum:compeuler} with the discretization given above, using the current stage solution as input. The source terms are then added as correction to the right-hand side.
The equation for the self-gravitational potential is not time-dependent but its state variable is depended upon by the source terms in Eq.~\eqref{eq:govnum:compeuler}, so the self-gravitational potential has to be solved for in each update of the hydrodynamics and also for intermediate steps. This is done by loading the coefficients for the mass density into a buffer on the C++ side of the implementation, solving for the self-gravitational potential and copying the coefficients for the self-gravitational potential back into buffers provided in the semi-discretization, these calling relations are further illustrated in Fig.~\ref{fig:algs:seq-rhs}.

The C++ solver for the self-gravitational potential is realized as a stand-alone application controlled through a wrapper class exchanging information with the Trixi.jl layer. From this class, the multigrid solver infrastructure is set up, defined on top of the usual management classes for enumerations of degrees of freedom of a DG field (in exact correspondence with Trixi.jl), labeled as \texttt{DoFHandler::distribute\_dofs(FE\_DGQ<dim>(degree\_dg))} and the continuous finite element field \texttt{DoFHandler::distribute\_dofs(FE\_Q<dim>(degree\_cg))}~\cite{arndt2021}, as well as parallel solution vectors, infrastructure for fast matrix-free evaluation~\cite{Kronbichler2012} and the multigrid levels~\cite{munch2023}. The library-driven design ensures that only those data structures needed for the coupled solver are put together, making the overall control code in the wrapper relatively lean.

For calling the deal.II-based multigrid solver, a clearly defined application binary interface to call into has to be provided, which C++ by default does not offer as there is no standard definition of the name mangling scheme. Hence, we provide a thin wrapping C library that just passes its arguments or pointers to its arguments through to either the multigrid solver or back to the calling Julia functions, so that the calling interface follows the standard C function calling convention. A handle pointing to an instance of the convenience class is then kept as parameter in the main semi-discretization. The calls into the appropriate methods of the C++ solver class we then wrap, again, in Julia functions for easier usage.

It should be noted that the default configurations of deal.II and Trixi.jl impose a different balancing on the mesh. While Trixi.jl only imposes the two-to-one balancing along edges or faces for the used mortar method, the deal.II-application imposes this balancing also over vertices. Hence, the balancing in Trixi.jl is adapted to the stronger balancing imposed in deal.II. Furthermore, the default ordering of unknowns in deal.II follows a level-order, which differs from the default Morton order used in {\ttfamily p4est} and hence vectors defined by Trixi.jl. For a consistent interpretation of unknowns, we impose a numbering strategy~\cite{Bangerth2007} to match the Morton order.

For enabling adaptive mesh refinement in the multigrid solver, the mesh refinement indicators computed in Trixi.jl need to be passed to deal.II, which then marks the respective mesh cells to be refined and performs the mesh adaptation following the general deal.II parallel mesh refinement procedures~\cite{Bangerth2011}.

\begin{figure}
    \centering
    \begin{tikzpicture}
        \draw[->] (-0.25,0) node[left] {coupled solver} -- (10.5,0);
        \draw[->] (-0.25,1.5) node[left] {hydro solver} -- (10.5,1.5);
        \draw[->] (-0.25,3.0) node[left] {gravity solver} -- (10.5,3.0);
        \filldraw[color=black,fill=white] (0,-0.25) -- (2,-0.25) -- (2,0.25) -- (0,.25) -- cycle;
        \filldraw[color=black,fill=white] (2,1.25) -- (3,1.25) -- (3,1.75) -- (2,1.75) -- cycle;
        \filldraw[color=black,fill=white] (3,-0.25) -- (4,-0.25) -- (4,0.25) -- (3,.25) -- cycle;
        \filldraw[color=black,fill=white] (4,2.75) -- (5,2.75) -- (5,3.25) -- (4,3.25) -- cycle;
        \filldraw[color=black,fill=white] (5,-0.25) -- (7,-0.25) -- (7,0.25) -- (5,.25) -- cycle;
        \filldraw[color=black,fill=white] (7,1.25) -- (9.5,1.25) -- (9.5,1.75) -- node[midway,below] {add source terms} (7,1.75) -- cycle;
        \filldraw[color=black,fill=white] (9.5,-0.25) -- (10,-0.25) -- (10,0.25) -- (9.5,.25) -- cycle;
        \draw[-Stealth] (0,-.75) -- node[midway,sloped,above]{rhs!} (0,-.25);
        \draw[-Stealth] (2,0.25) -- node[midway,sloped,above]{rhs!} (2,1.25);
        \draw[-Stealth] (3,1.25) -- (3,0.25);
        \draw[->] (1,0.35) -- (1,1.40);
        \draw[->,dashed] (1,1.60) -- node[midway,sloped,above]{store \(\rho\)}(1,2.90);
        \draw[-Stealth] (4,0.25) -- node[midway,sloped,above]{solve}(4,2.75);
        \draw[-Stealth] (5,2.75) -- (5,0.25);
        \draw[->] (6,0.35) -- (6,1.4);
        \draw[color=white] (6,1.6) -- node[color=black,midway,sloped,above]{load \(\nabla\phi\)} (6,2.9);
        \draw[->,dashed] (6,2.9) -- (6,1.6);
        \draw[-Stealth] (7,0.25) -- (7,1.25);
        \draw[-Stealth] (9.5,1.25) -- (9.5,0.25);
        \draw[-Stealth] (10,-.25) -- (10,-.75);

    \end{tikzpicture}
    \caption{Sequence diagram for the evaluation of the semi-discretization for the coupled system}
    \label{fig:algs:seq-rhs}
\end{figure}
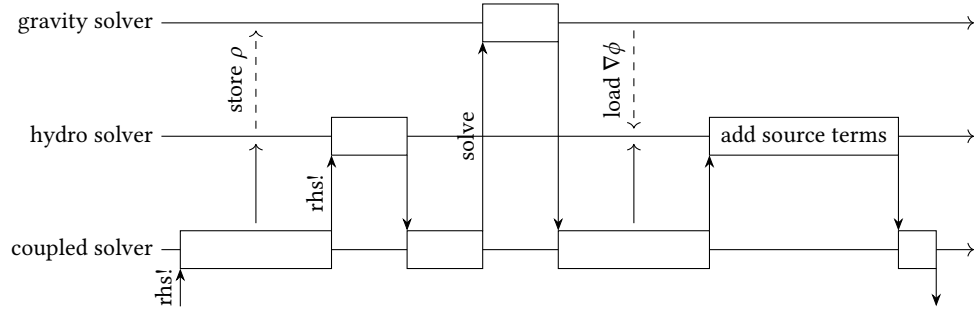

\begin{figure}
    \centering
    \begin{tikzpicture}
        \draw[->] (-0.25,0) node[left] {AMR callback} -- (10.5,0);
        \draw[->] (-0.25,1.5) node[left] {mesh} -- (10.5,1.5);
        \draw[->] (-0.25,3.0) node[left] {hydro solver} -- (10.5,3.0);
        \draw[->] (-0.25,4.5) node[left] {gravity solver} -- (10.5,4.5);
        \filldraw[color=black,fill=white] (0,-0.25) -- (4.0,-0.25) -- (4.0,0.25) -- (0,.25) -- cycle;
        \filldraw[color=black,fill=white] (4.0,1.25) -- (4.5,1.25) -- (4.5,1.75) -- (4.0,1.75) -- cycle;
        \filldraw[color=black,fill=white] (4.5,-0.25) -- (6.5,-0.25) -- (6.5,0.25) -- (4.5,.25) -- cycle;
        \filldraw[color=black,fill=white] (6.5,3.25) -- (7,3.25) -- (7,2.75) -- (6.5,2.75) -- cycle;
        \filldraw[color=black,fill=white] (7,-.25) -- (9,-.25) -- (9,.25) -- (7,.25) -- cycle;
        \filldraw[color=black,fill=white] (9,4.25) -- (9,4.75) -- (9.5,4.75) -- (9.5,4.25) -- cycle;
        \filldraw[color=black,fill=white] (9.5,-0.25) -- (10,-0.25) -- (10,0.25) -- (9.5,.25) -- cycle;
        \draw[-Stealth] (0,-.75) -- (0,-.25);
        \draw[loosely dotted,color=white] (.5,1.6) -- node[color=black,above,midway,sloped]{\(\rho\)} (.5,2.9);
        \draw[dashed,->] (.5,2.9) -- (.5,.35);
        \draw[->,dashed] (3.25,0.25) --node[midway,sloped,above]{inds} (3.25,1.4);
        \draw[-Stealth] (4,0.25) --node[midway,sloped,above]{adapt} (4,1.25);
        \draw[-Stealth] (4.5,1.25) -- (4.5,0.25);
        \draw[->,dashed] (5.5,0.25) --node[midway,sloped,above]{\;adapted cell list} (5.5,2.9);
        \draw[color=white] (6.5,1.8) --node[color=black,midway,sloped,above] {adapt} (6.5,2.7);
        \draw[-Stealth] (6.5,0.25) -- (6.5,2.75);
        \draw[-Stealth] (7.0,2.75) -- (7.0,0.25);
        \draw[color=white] (8.0,1.8) --node[color=black,midway,sloped,above] {indicators} (8.0,2.7);
        \draw[->,dashed] (8.0,0.25) -- (8.0,4.4);
        \draw[color=white] (9.0,1.8) --node[color=black,midway,sloped,above] {adapt} (9.0,2.7);
        \draw[-Stealth] (9.0,0.25) -- (9.0,4.25);
        \draw[-Stealth] (9.5,4.25) -- (9.5,0.25);
        \draw[-Stealth] (10,-.25) -- (10,-.75);
        \node at (2.0,0.0) {compute indicators};

    \end{tikzpicture}
    \caption{Sequence diagram for the adaptive mesh refinement process in the respective callback.}
    \label{fig:algs:seq-amr}
\end{figure}
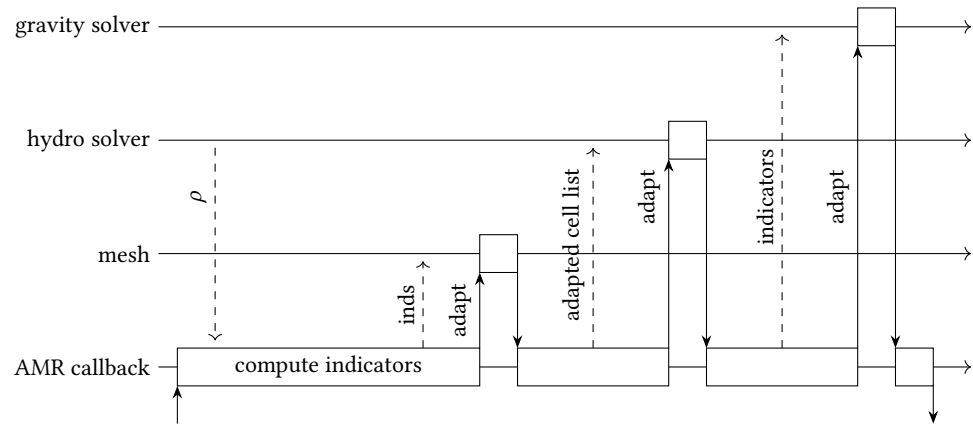

\begin{figure}
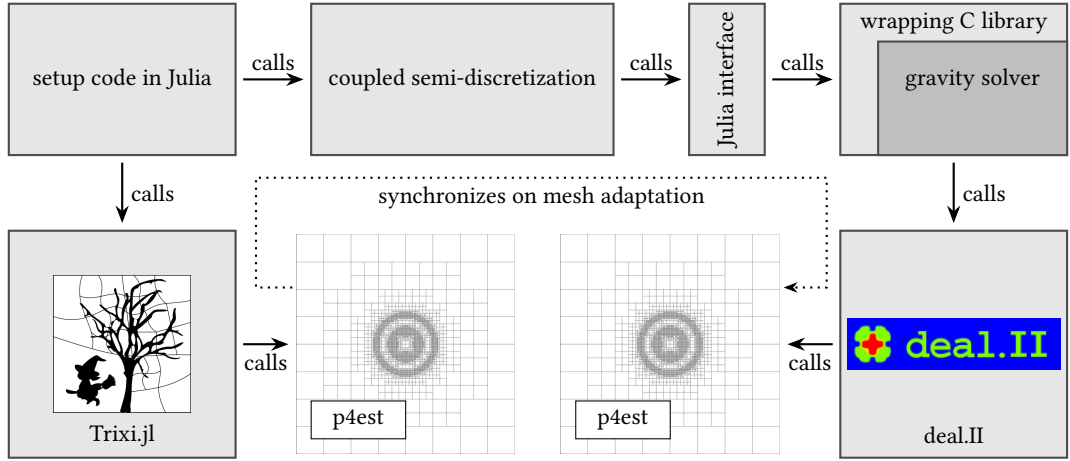

    \centering
    \begin{tikzpicture}
      \filldraw[thick,color=white!30!black,fill=white!90!black] (0,0)--(3,0)--(3,3)--(0,3)--cycle;
      \node at (1.5,0.3) {Trixi.jl};
      \node at (1.5,1.5) {\includegraphics[scale=0.12]{trixi.png}};
      \node at (5.25,1.5) {\includegraphics[scale=0.8]{p4est_mesh.png}};
      \filldraw[fill=white] (4.0,0.25) -- (5.25,0.25) -- (5.25,0.75) -- (4.0,0.75) -- cycle;
      \node at (4.625,0.50) {p4est};
      \node at (8.75,1.5) {\includegraphics[scale=0.8]{p4est_mesh.png}};
      \filldraw[fill=white] (7.5,0.25) -- (8.75,0.25) -- (8.75,0.75) -- (7.5,0.75) -- cycle;
      \node at (8.125,0.50) {p4est};
      \filldraw[thick,color=white!30!black,fill=white!90!black] (11,0)--(14,0)--(14,3)--(11,3)--cycle;
      \node at (12.5,0.3) {deal.II};
      \node at (12.5,1.5) {\includegraphics[scale=0.4]{dealii.png}};

      \draw[thick,-Stealth] (3.1,1.5) -- node[midway,below]{calls} (3.7,1.5);
      \draw[thick,Stealth-] (10.3,1.5) -- node[midway,below]{calls} (10.9,1.5);
      \filldraw[thick,color=white!30!black,fill=white!90!black] (11,4)--(14,4)--(14,6)--(11,6)--cycle;
      \node at (12.5,5.75) {wrapping C library};
        \filldraw[thick,color=black!70!white,fill=black!25!white] (11.50,4)--(14,4)--(14,5.5)--(11.50,5.5)--cycle;
        \node at (12.75,5.00) {gravity solver};
      \draw[thick,Stealth-] (12.5,3.1) -- node[midway,right]{calls} (12.5,3.9);
      \filldraw[thick,color=white!30!black,fill=white!90!black] (0,4)--(3,4)--(3,6)--(0,6)--cycle;
      \node at (1.5,5) {setup code in Julia};
      \filldraw[thick,color=white!30!black,fill=white!90!black] (4,4)--(8,4)--(8,6)--(4,6)--cycle;
      \node at (6,5) {coupled semi-discretization};
      \draw[thick,-Stealth] (8.1,5.0) -- node[midway,above] {calls} (8.9,5.0);
      \filldraw[thick,color=white!30!black,fill=white!90!black] (9,4)--(10,4)--(10,6)--(9,6)--cycle;
      \draw[color=white!90!black] (9.5,4.1)-- node[color=black,sloped,midway]{Julia interface} (9.5,5.9);
      \draw[thick,-Stealth] (3.1,5.0) -- node[midway,above] {calls} (3.9,5.0);
      \draw[thick,Stealth-] (1.5,3.1) -- node[midway,right] {calls} (1.5,3.9);
      \draw[thick,-Stealth] (10.1,5.0) -- node[above,midway] {calls} (10.9,5.0);
      \draw[thick,arrows={-Stealth[scale=0.8]},dotted] (3.75,2.25) -- (3.27,2.25) -- (3.27,3.7) -- (10.82,3.7) node[below,midway] {synchronizes on mesh adaptation} -- (10.82,2.25) -- (10.25,2.25);
   \end{tikzpicture}
    \Description{}
    \caption{A schematic illustration of the coupling and calling relations during a numerical simulation between the components comprising the coupled solver.}
    \label{fig:algs:wrapping}
\end{figure}

\subsection{Reproducible simulations}

As the main control flow in the numerical simulation is on the side of the implementation in Julia, we can directly use Julia's project and dependency description in the {\ttfamily Project.toml} and {\ttfamily Manifest.toml} files to provide machine-readable files to retrieve the settings used to produce results in a machine-independent fashion. This is possible as these files contain the transitive dependencies of the project and also include which versions were used to obtain the results in this work, along with possible compatibility bounds of the dependencies. However, as this only covers Julia packages, we provide further packages containing a build of the deal.II library and of the coupled multigrid solver. The deal.II package with a choice of common dependencies is now available in the global registry of the Julia package manager. For building these packages we can, however, not just refer to the builds already provided by the developers of deal.II as we need to properly set up the dependencies of deal.II to other Julia packages implementing the desired dependencies.

The environment for building Julia packages is provided by the BinaryBuilder.jl package, which provides a sand-boxed, minimally outfitted system running a Linux distribution with a fixed tree structure. The artifacts contained in the packages are cross-compiled for each of the target platforms and then packaged as collections of text files and binary blobs for the different platforms. These packages can be deployed and registered globally, enabling each Julia user to add the package to their project. Alternatively, packages can be deployed to an arbitrary accessible git repository which can then be added as dependency to projects without taking up a space in the namespace of the global registry. 

We provide a reproducibility repository in \cite{ehlert2026multisolverRepro}, which contains a Julia project depending on the coupling package containing the implementation of the multigrid solver and its dependency to the deal.II package as added to the global Julia package registry. From this project, we can then call into the artifact of the compiled library inside the package containing the coupling code to call the proper functions we implemented therein to initialize, compute a solution, adapt the mesh or release the memory allocated by the wrapper for the multigrid solver.

As we can build Julia packages for deal.II, and subsequently also for our multigrid solver contained in its wrapping C library for some of the more common platforms supported by Julia, this does not only yield reproducible runs of the coupled solvers, but also portability of our coupled solver across the more common operating systems, which include MacOS and Linux distributions on common modern hardware architectures. As modern Microsoft Windows operating systems already offer a Linux compatibility layer, the Windows Subsystem for Linux, the proposed solver is portable to each of these operating systems, obviating a native build there.

\section{Numerical Results} \label{sc:numres}


In this section, we perform a set of numerical experiments with our code. First, we verify our code against a convergence test to ensure its correctness. Then, we apply it to simulate the Jeans instability to verify our solver against a physical problem for which there are analytical energy profiles to compare against. Afterwards, we simulate a cylindrical Sedov blast wave in two spatial dimensions to show interoperability with enabled \(h\)-adaptivity and shock-capturing in the coupled solver as well as checking the computational performance and scaling of our code by applying it to simulate spherical Sedov blast waves in three spatial dimensions.

A reproducibility repository \cite{ehlert2026multisolverRepro} contains the implementation of our coupled solver, scripts to reproduce the results presented in this section as well as the actual results. The results in this section have been tested to be reproducible using Julia version 1.12.4 on Linux and MacOS running on either {\ttfamily aarch64} or {\ttfamily x86\_64} platforms.

\subsection{Convergence test} \label{sc:numres:conv}

For verifying the high-order numerical accuracy of the implemented solvers, we test the coupled solver with equations derived using the method of manufactured solutions. We consider the system of Eq.~\eqref{eq:govnum:compeuler} and Eq.~\eqref{eq:govnum:selfgrav} in two dimensions on a domain \(\Omega=[0,2]^2\) along with periodic boundary conditions and \(\gamma=2\). We specify a  solution of the form
\begin{equation}
    \rho = 2 + \frac1{10}\sin(\pi(x+y-t)),\quad v_1=v_2=1,\quad p = \frac1\pi\rho^2, \quad \phi = -\frac2\pi(\rho -2). \label{eq:numres:conv:solution}
\end{equation}
For the derivatives of the state variables we get the following relations by the governing equations
\begin{equation}
    \rho_x = \rho_y = -\rho_t = \frac\pi{10}\cos(\pi(x+y-t)), \quad \phi_x = \phi_y = -\frac15(\cos(\pi(x+y-t)), \label{eq:numres:conv:derivative_relations}
\end{equation}
which are also periodic in the domain. Furthermore, there holds
\begin{equation}
    -\Delta \phi = -(\phi_{xx} + \phi_{yy}) = -4\pi(\rho - 2) = -4\pi\rho + 8\pi. \label{eq:numres:conv:laplace}
\end{equation}
This setup solves the Poisson equation for the self-gravitational potential with gravity constant \(G=1\) and constant in time residual \(8\pi\). For the self-gravitational equations the residuals in total can then be determined to be
\begin{equation}
    \ubf_t + \nabla\cdot\fbf(\ubf) = \sbf(\ubf) + \begin{bmatrix}\frac\pi{10}\cos(\pi\,c) \\[2pt] \frac\pi{10}\cos(\pi\,c) \\[2pt] \frac\pi{10}\cos(\pi\,c) \\[2pt] \frac\pi{10}\cos(\pi\,c)(1+\frac2\pi(2+\frac1{10}\sin(\pi\,c)))\end{bmatrix},\label{eq:numres:conv:residuals}
\end{equation}
where \(c = x+y-t\).

This test case is run until a final time of \(t_\mathrm{end}=1\) on a sequence of uniformly refined Cartesian meshes with a bisection strategy. We investigate the expected order of convergence of the approximation in the DG space with two different polynomial orders for the DGSEM of three and four. The discrete \(L^2\) errors of the state variables are computed using LGL-quadrature at the end of the test runs over the entire domain \(\Omega\).

The results of this test case can be found in Tab.~\ref{tab:numres:conv:eoc}, including the deduced experimental order of convergence. We obtain values close to the predicted order of convergence \(p+1\) for a DGSEM with a piecewise polynomial approximation space of order \(p\). This verifies the coupled solver.

\begin{table}
    \centering
    polynomial order of \(3\)
    
    \begin{tabular}{ccccc}
        \toprule
        mesh elements & \(L^2(\rho)\) & \(L^2(\rho\,v_1)\) & \(L^2(\rho\,v_2)\) & \(L^2(E)\) \\
        \midrule
        \(2^2\) & \(5.84\cdot10^{-3}\) & \(6.72\cdot10^{-3}\) & \(6.72\cdot10^{-3}\) & \(1.34\cdot10^{-2}\) \\[2pt]
        \(4^2\) & \(2.97\cdot10^{-4}\) & \(3.65\cdot10^{-4}\) & \(3.65\cdot10^{-4}\) & \(6.71\cdot10^{-4}\) \\[2pt]
        \(8^2\) & \(1.97\cdot10^{-5}\) & \(2.33\cdot10^{-5}\) & \(2.33\cdot10^{-5}\) & \(4.28\cdot10^{-5}\)\\[2pt]
        \(16^2\) & \(1.03\cdot10^{-6}\) & \(1.41\cdot10^{-6}\) & \(1.41\cdot10^{-6}\) & \(2.65\cdot10^{-6}\)\\[2pt]
        \(32^2\) & \(6.38\cdot10^{-8}\) & \(8.80\cdot10^{-8}\) & \(8.80\cdot10^{-8}\) & \(1.65\cdot10^{-7}\)\\
        \midrule
        avg. EOC & \(4.12\) & \(4.06\) & \(4.06\) & \(4.08\) \\
        \bottomrule 
    \end{tabular}
    
    \vspace{\baselineskip}
    polynomial order of \(4\)
    
    \begin{tabular}{ccccc}
        \toprule
        mesh elements & \(L^2(\rho)\) & \(L^2(\rho\,v_1)\) & \(L^2(\rho\,v_2)\) & \(L^2(E)\) \\
        \midrule
        \(2^2\) & \(6.59\cdot10^{-4}\) & \(8.19\cdot10^{-4}\) & \(8.19\cdot10^{-4}\) & \(1.51\cdot10^{-3}\) \\[2pt]
        \(4^2\) & \(2.21\cdot10^{-5}\) & \(2.65\cdot10^{-5}\) & \(2.65\cdot10^{-5}\) & \(4.76\cdot10^{-5}\) \\[2pt]
        \(8^2\) & \(7.85\cdot10^{-7}\) & \(9.01\cdot10^{-7}\) & \(9.01\cdot10^{-7}\) & \(1.66\cdot10^{-6}\) \\[2pt]
        \(16^2\) & \(1.90\cdot10^{-8}\) & \(2.59\cdot10^{-8}\) & \(2.59\cdot10^{-8}\) & \(4.70\cdot10^{-8}\) \\[2pt]
        \(32^2\) & \(6.18\cdot10^{-10}\) & \(8.19\cdot10^{-10}\) & \(8.19\cdot10^{-10}\) & \(1.48\cdot10^{-9}\) \\
        \midrule
        avg. EOC & \(5.01\) & \(4.98\) & \(4.98\) & \(4.99\) \\
        \bottomrule 
    \end{tabular}
    \Description{Two tables showing the L2 errors of the numerical simulation of the problem with the manufactured simulation in the current subsection for polynomial orders of three (3) and four (4). The average experimental order of convergence is consistent with what we expect from the DGSEM that it is one order higher than the order of the polynomials used to approximate the solution space, with a bit of wiggling.}
    \caption{Results of the convergence test for the manufactured solution in subsection~\ref{sc:numres:conv}.}
    \label{tab:numres:conv:eoc}
\end{table}

\subsection{Jeans instability} \label{sc:numres:jeans}

To further validate the model and our solver with an actual physical test case, we consider the Jeans instability in two dimensions, which describes an instability in a self-gravitating, thermally supported interstellar cloud which was first described by Jeans \cite{jeans1902}. For this case, the linear instability mode with analytic energy profiles \cite{flashuserguide, derigs2016} is especially useful for testing our coupled solver.

For the setup of this problem we follow \cite{derigs2016}, see there for further details. Along with periodic boundary conditions, the initial condition is given by
\begin{align}
    \begin{cases}
    \rho(\xbf) = \rho_0\,(1+\delta\cos(\kbf\cdot\xbf)) & \forall x\in\Omega,
    \\[-2pt]
    p(\xbf) = p_0\,(1+\gamma\,\delta\cos(\kbf\cdot\xbf)) & \forall x\in\Omega,
    \\[-2pt]
    \vbf(\xbf) = 0 & \forall x\in\Omega.
    \end{cases}
\end{align}
Here \(\delta\ll1\) is the amplitude of the perturbation. Whether the perturbation is stable is determined by the relationship between the wavenumber \(k = |\kbf|=\|\kbf\|_2\) and the Jeans wavenumber \(k_J\), given by
\begin{align}
    k_J = \frac{\sqrt{4\pi\,G\,\rho_0}}{c_0},
\end{align}
where \(c_0 = \sqrt{\gamma\,p_0/\rho_0}\) is the unperturbed adiabatic speed of sound \cite{chandrasekhar1961}. For \(k_J < k\) the dynamics are stable and the perturbation oscillates with a frequency of
\begin{align}
    \omega = \sqrt{c_0^2\,k^2 - 4\pi\,G\,\rho_0 } = \sqrt{c_0^2\,(k^2 - k_J^2)},
\end{align}
while for \(k_J > k\) the term under the square root would get negative and yield an imaginary \(\omega\) which would correspond to unstable dynamics. The heat capacity ratio is given by \(\gamma=\frac53\). As in \cite{derigs2016} we use \(p_0 = 1.5\cdot10^7\;[\mathrm{dyn}\,\mathrm{cm}^{-2}]\) and \(\rho_0=1.5\cdot10^7\;[\mathrm{g}\,\mathrm{cm}^{-3}]\) which yields a Jeans wavenumber of \(k_J=2.747\;[\mathrm{cm}^{-1}]\).

As the domain is periodic, we can choose a domain size and resolution to avoid any unstable frequencies getting amplified by round-off errors due to floating-point computations. This is fulfilled when \(k_J\) is smaller than the smallest resolvable wave number on the domain. We choose the domain \(\Omega\) to fulfill \(\Omega \subset [0,L]^2\) where \(L = \frac\pi{k_J} = \frac12\sqrt{\frac{\pi\,\gamma\,p_0}{G\,\rho_0^2}}\), with the restriction that a full period of the initial condition is sampled on the domain.

\begin{figure}
    \includegraphics[scale=0.5]{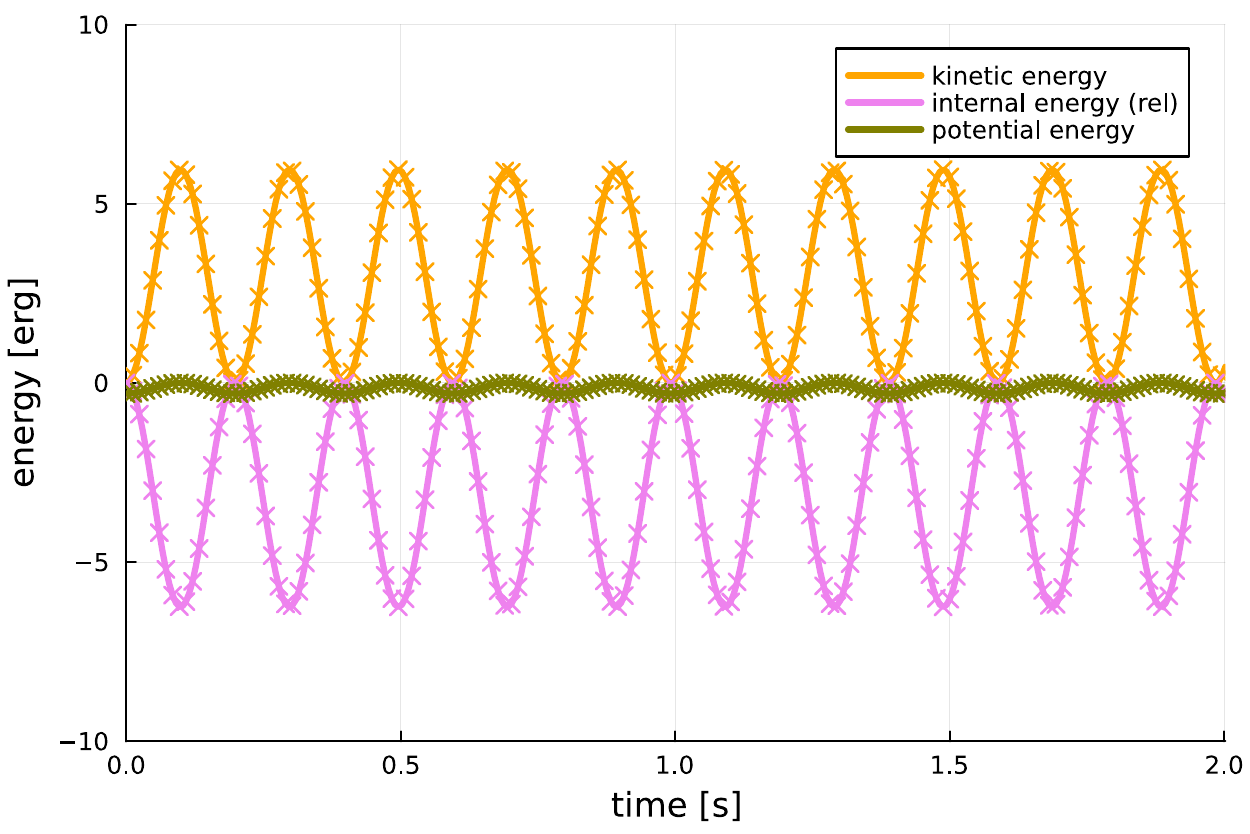}
    \caption{The continuous lines show the analytic values of the energies, the cross marks show the energies at every second time step during the numerical simulation of the Jeans instability with wave vector \(k=(\pi, 0)\)}
    \Description{Figures of numerical simulations of the Jeans instability for two seconds using the Trixi.jl--deal.II coupling presented in this work. Plotted are the kinetic, relative internal and potential energy of every second time step of the numerical simulation of the Jeans instability along with their analytical values which align perfectly.}
    \label{fig:numres:jeans:comparisons}
\end{figure}

The analytical expressions for the kinetic, relative internal and potential energies are given by
\begin{align}
    E_\mathrm{kin}(t) &= \frac{\rho_0\,\delta^2\,|\omega|^2\,L^2}{8k^2}\left[1-\cos(2\omega\,t)\right], \label{eq:numres:ekin}
    \\[2pt]
    E_\mathrm{int}(t) - E_\mathrm{int}(0) &= -\frac18\,\rho_0\,c_0^2\,\delta^2\,L^2\left[1-\cos(2\omega\,t)\right], \label{eq:numres:eint}
    \\[2pt]
    E_\mathrm{pot}(t) &= -\frac{\pi\,G\,\rho_0^2\,\delta^2\,L^2}{2k^2}\left[1+\cos(2\omega\,t)\right]. \label{eq:numres:epot}
\end{align}
The results of the numerical simulation are evaluated by comparing the kinetic, internal and potential energies of the dynamical system against these analytical energies. The kinetic, internal and potential energies are computed by
\begin{align}
    E_\mathrm{kin} = \int \frac\rho2\left(v_1^2 + v_2^2\right)\;d\Omega, \quad E_\mathrm{int} = \int \frac{p}{\gamma - 1}\;d\Omega, \quad E_\mathrm{pot} = \frac12\int \rho\,\phi\;d\Omega.
\end{align}
These integrals are determined numerically by LGL quadrature, as the dynamical system is already discretized for the DGSEM using LGL nodes. For the comparison against analytical values see Fig.~\ref{fig:numres:jeans:comparisons}, where one can notice that the energies computed in the numerical simulation align near perfectly with the analytical energies.

\subsection{Sedov blast wave} \label{sc:numres:sedov}

As a next test case, we consider cylindrical and spherical Sedov blast waves  \cite{sedov2018}, which can be understood in a self-gravitational setting to be strongly simplified models of super novae. In this test case, we use the shock-capturing and \(h\)-adaptivity features provided by Trixi.jl embedded into the proposed strategy of subsection~\ref{sc:algs:technicals}.

The Sedov blast wave is the result of a strong explosion, the energy of which is deposited into a single point in a medium of ambient pressure \(p_\mathrm{am}\) and density \(\rho_\mathrm{am}\). For practical simulations, this model is not useful because the energetic area typically can not be resolved by the mesh, which is why we  follow \cite{schlottkelakemper2021,flash} and convert the energy of the explosion to a pressure which is deposited in a small ball of radius \(r_\mathrm{ini}\) centered on the origin. The radius of this ball is chosen to be four times the element size at the finest resolution of the mesh to avoid imprinting of the mesh on the initial values. The energy is converted to a pressure value by
\begin{align}
    p_\mathrm{ini}(r)\;[\mathrm{dyn\;cm^{-2}}] = \begin{cases} 
    \frac{3(\gamma-1)\,E}{(d+1)\pi\,r_\mathrm{ini}^d}, & r \leq r_\mathrm{ini} \\
    p_\mathrm{am}, & r > r_\mathrm{ini},
    \end{cases} \label{eq:numres:sedov:inipres}
\end{align}
where \(d\) is the dimension of the space in which we simulate.

We give the parameters of the Sedov blast wave in CGS units, where \(p_\mathrm{am} = 10^{-5}\;[\mathrm{dyn\;cm^{-2}}]\), \(E = 1\;[\mathrm{erg}]\), \(v_1=v_2=0\;[\mathrm{cm\;s^{-1}}]\) and the heat capacity ratio \(\gamma = 1.4\). The computational domain is \(\Omega_2=[-4,4]^2\;[\mathrm{cm^2}]\), which is discretized with an adaptive mesh that is allowed a finest resolution of \(2^{-5}\;[\mathrm{cm}]\), i.e.~a maximum refinement level of eight. Hence, \(r_\mathrm{ini}\) takes the value of \(2^{-3}\;[\mathrm{cm}]\) in equation~\eqref{eq:numres:sedov:inipres}. The gravitational constant \(G\) has the value \(6.674\cdot10^{-8}\;[\mathrm{cm^3}\;g^{-1}\;s^{-2}]\). As the gravitational potential is usually assumed to vanish at large distances away from a localized region containing a non-vanishing gravitational potential \cite{flashuserguide}, we localize the ambient density of the blast wave inside a ball with radius \(r_\rho=1\;[\mathrm{cm}]\), such that
\begin{align}
    \rho_\mathrm{am}(r)\;[\mathrm{g\,cm^{-3}}] = \begin{cases}
        1, & r \leq r_\rho, \\ 10^{-5}, & r > r_\rho.
    \end{cases}
\end{align}
The transition between inner and ambient state for both density and pressure is smoothed by a logistic function \(\sigma_{k,r}\) with steepness \(k=150\;[\mathrm{cm^{-1}}]\), centered on  the respective radii \(r_\rho\) or \(r_\mathrm{ini}\). Ultimately, the initial values for a cylindrical Sedov blast wave are given by:
\begin{align}
    \begin{cases}
        \rho(x,0) = \sigma_{150,r_\rho}(\rho_\mathrm{am}(\|x\|)) & \forall x\in\Omega_2,
        \\[-2pt]
        p(x,0) = \sigma_{150,r_\mathrm{ini}}(p_\mathrm{ini}(\|x\|)) & \forall x\in\Omega_2,
        \\[-2pt]
        v_1(x,0) = v_2(x,0) = 0 & \forall x\in\Omega_2.
    \end{cases} \label{eq:numres:inivalssedov2d}
\end{align}
For a spherical Sedov blast wave, we consider \(\Omega_3 = [-2,2]^3\) as computational domain, which is discretized with an adaptive mesh with the same finest resolution as in the two dimensional case and hence with a maximum refinement level of seven, to reduce the main memory requirements of the simulation. The initial values considered in three dimensions are given by:
\begin{align}
    \begin{cases}
        \rho(x,0) = \sigma_{150,r_\rho}(\rho_\mathrm{am}(\|x\|)) & \forall x\in\Omega_3,
        \\[-2pt]
        p(x,0) = \sigma_{150,r_\mathrm{ini}}(p_\mathrm{ini}(\|x\|)) & \forall x\in\Omega_3,
        \\[-2pt]
        v_1(x,0) = v_2(x,0) = v_3(x,0) = 0 & \forall x\in\Omega_3.
    \end{cases} \label{eq:numres:inivalssedov3d}
\end{align}
We approximate the hydrodynamic and gravitational equations in both two and three dimensions with a polynomial degree of three and perform a numerical simulation until a final time \(t_{\mathrm{end},2}=2.2\;[\mathrm{s}]\) in the two-dimensional case and until a final time \(t_{\mathrm{end},3}=\frac12\;[\mathrm{s}]\) in the three-dimensional case. Time steps are selected according to a CFL-condition of \(1.0\). For the compressible Euler solver, we use the split-form DGSEM feature of Trixi.jl along with the shock capturing method described in subsection~\ref{sc:govnum:dgsem}. As indicator for mesh refinement the same indicator function as for shock capturing is evaluated, marking each cell containing a shock for refinement.

For a reference profile of the density to compare against, we simulated the cylindrical Sedov explosion on a uniform mesh refined to the finest resolution that was allowed in the \(h\)-adaptive simulation run just before.

\begin{figure}
    \centering
    \includegraphics[width=0.49\linewidth]{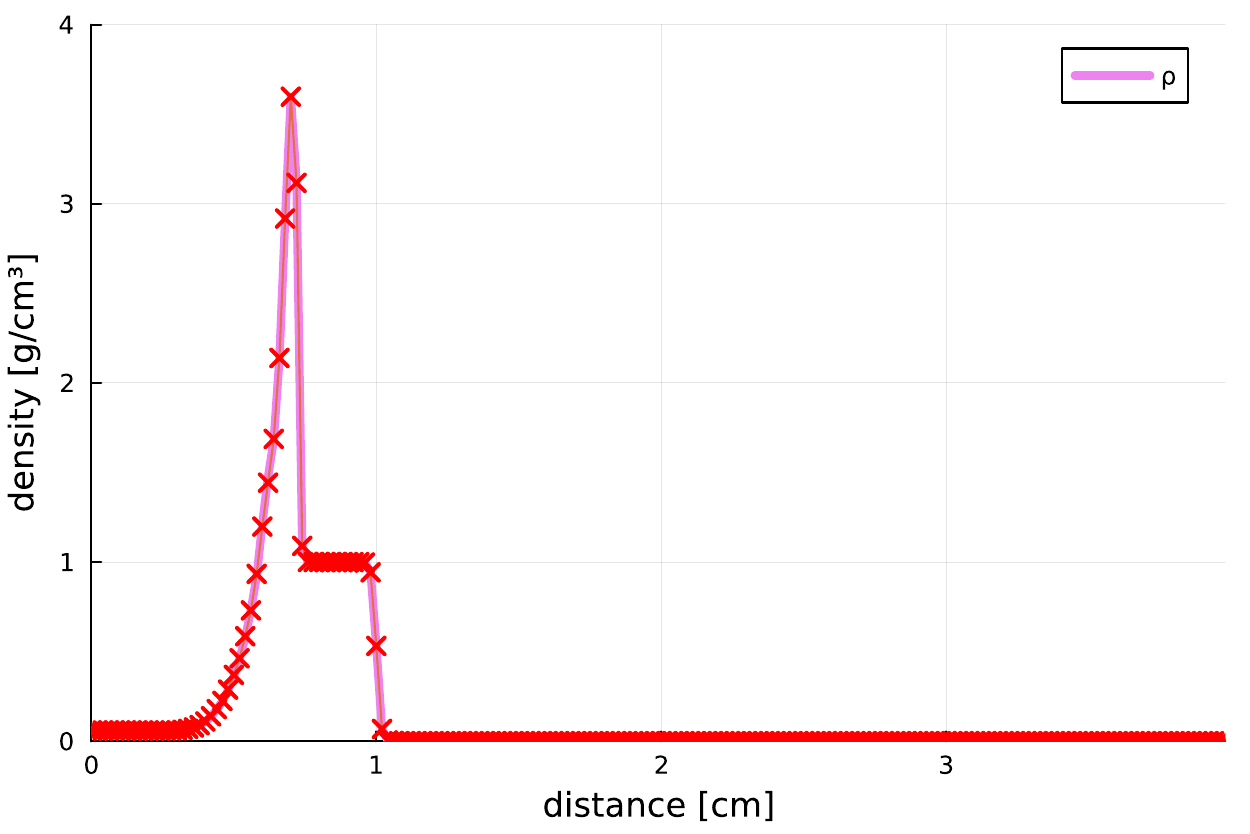}
    \includegraphics[width=0.49\linewidth]{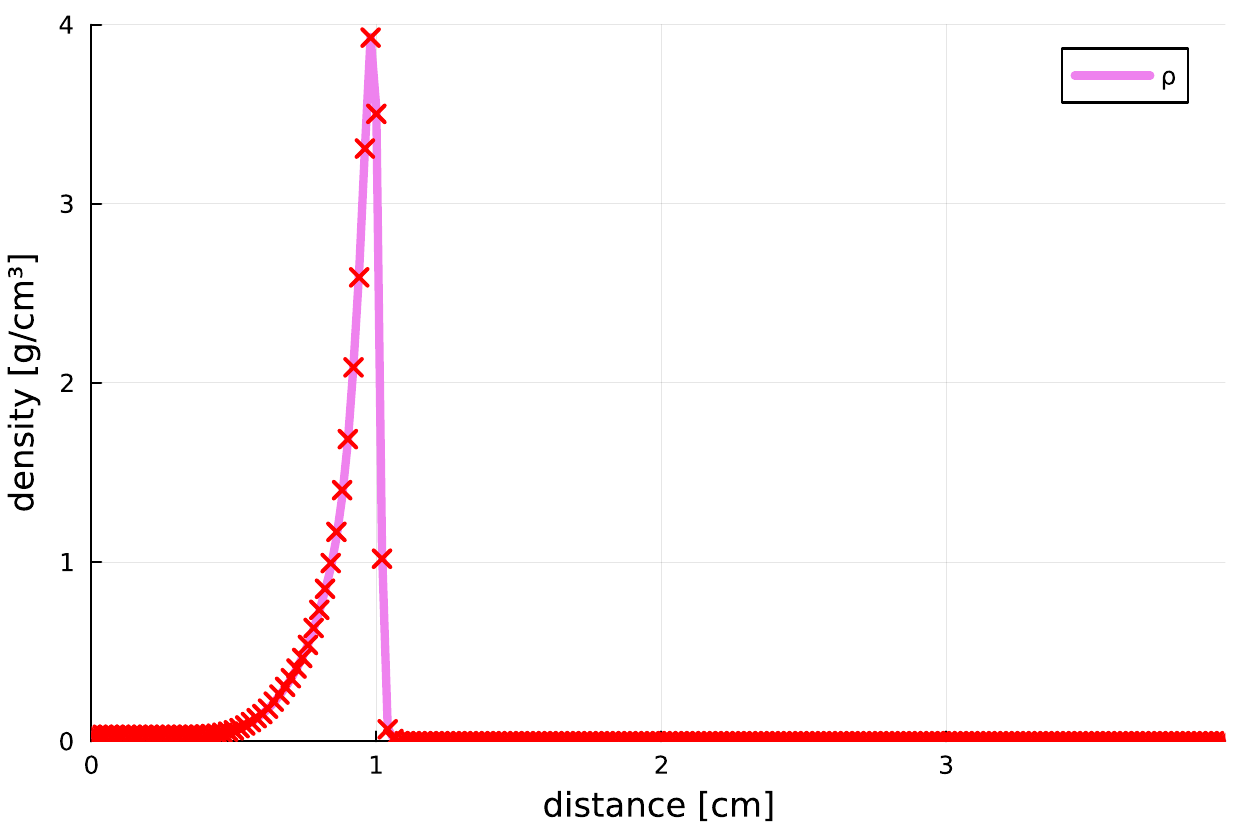}

    \vspace{4mm}

    \setlength{\sedovw}{0.49\linewidth}
    \sedovsplit{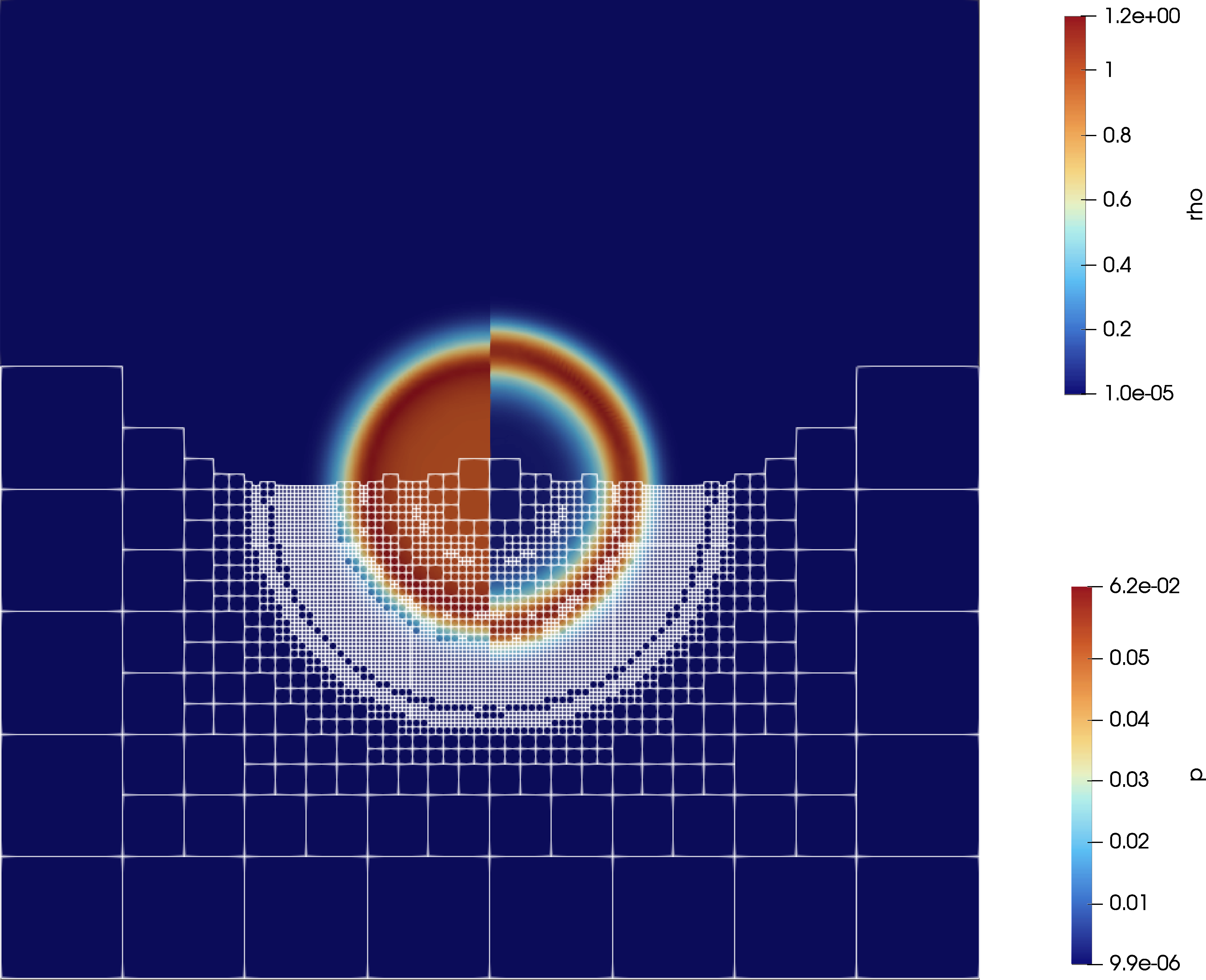}
    \sedovsplit{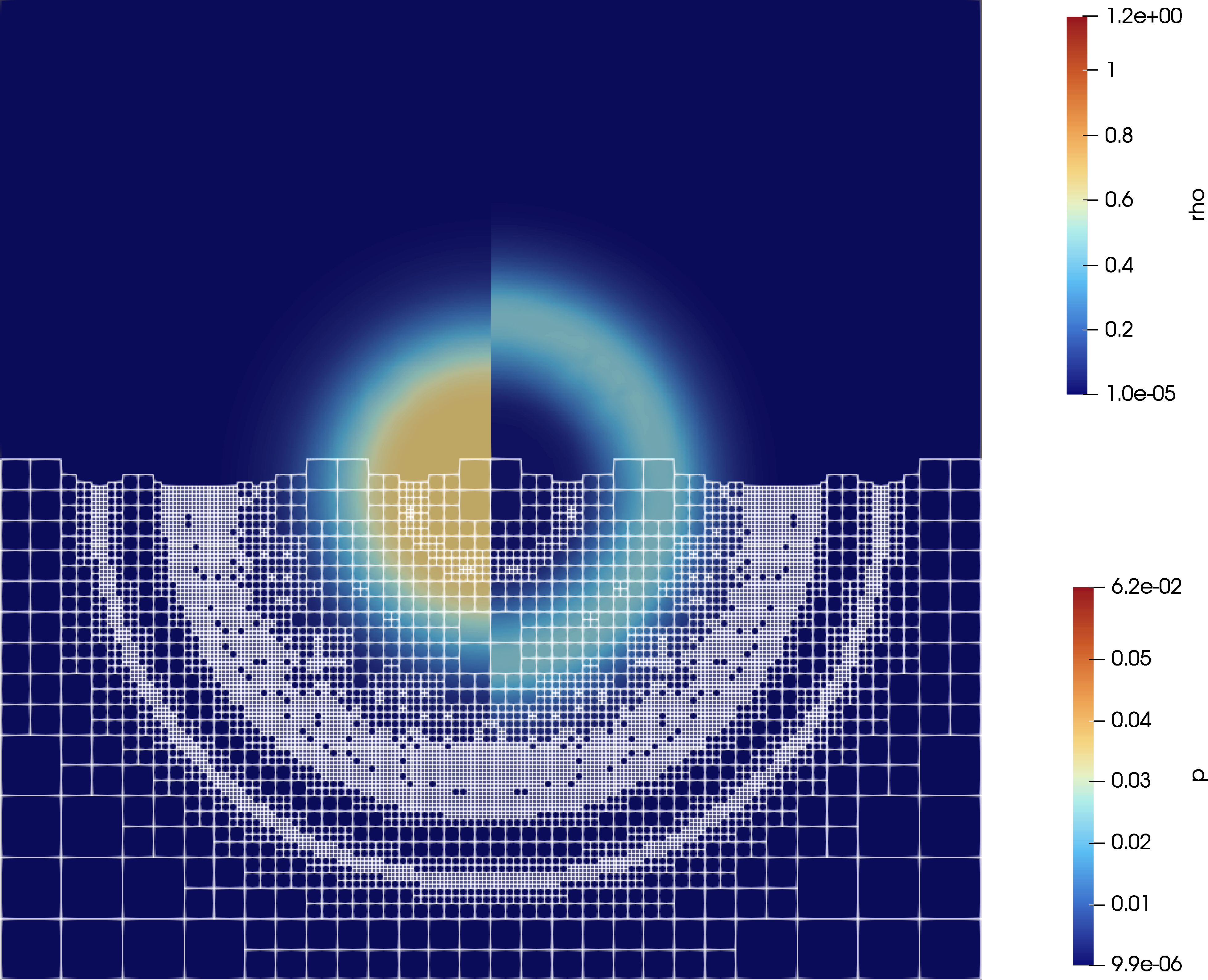}
    \caption{Top row: Radial density profiles of the cylindrical Sedov blast wave sampled along the positive x-axis after simulating half a second (left) and one second (right).
    Bottom row: Plots showing pressure (left half) and density (right half) of a cylindrical Sedov blast wave after \(1.5\,[\mathrm{s}]\) (left) and \(2.2\,[\mathrm{s}]\) (right), with an overlay of the adaptive mesh used in the simulation on the lower half.
    \label{fig:numres:sedov:profile}}
    \Description{In this figure, we show the radial density profiles of the numerical simulations of the cylindrical Sedov blast wave. Plotted are the values of a simulation with adaptive mesh refinement and a simulation without which, however, was computed on a mesh uniformly refined to the finest resolution allowed in the adaptive simulation. Visually, there are no recognizable differences between the density values of the adaptive and the uniformly refined simulations.}
\end{figure}

In Fig.~\ref{fig:numres:sedov:profile} we can see that the \(h\)-adaptive simulation run agrees near perfectly with the output of the uniformly refined mesh, as is the expected result. There is also a noticeable performance improvement, which we, however, will remark more on in the next subsection regarding computational performance.

\subsection{Computational performance} \label{sc:numres:perf}

To determine the performance and the scaling of our code, we perform a benchmark of the coupled solver by numerically simulating three dimensional Sedov blast waves. The timings of these benchmarks are taken on a machine running Ubuntu Noble on an AMD Ryzen Threadripper 3990X 64-Core Processor with \(251 \mathrm{GiB}\) of RAM and using MPICH as MPI implementation.
We consider the initial value problem for Eqs.~\eqref{eq:govnum:compeuler} and \eqref{eq:govnum:selfgrav} given by Eq.~\eqref{eq:numres:inivalssedov3d} on the domain \(\Omega_3=[-2,2]^3\) along with the same setup for the numerical simulation described in subsection~\ref{sc:numres:sedov}. For the test cases in which uniformly refined meshes are used, the domain is \([-1,1]^3\) to further reduce main memory requirements, also we then set \(r_\rho=\frac12\) in the initial values in Eq.~\eqref{eq:numres:inivalssedov3d}. Each of the numerical simulations is performed until the final time of \(t_\mathrm{end} = \frac12\;[\mathrm{s}]\) is reached, where the time step size is selected according to a CFL-condition of \(1.0\).

We investigate the scaling of the coupled solver on a mesh uniformly refined to level six, the total number of degrees of freedom in this numerical simulation is then \(2^{6\cdot3}\cdot 4^3 = 2^{24} \approx 16.8\,\text{M}\), as there are four nodes in each dimension for each of the cells of the uniformly discretized domain. Furthermore, we investigate the scaling of the coupled solver with enabled adaptive mesh refinement after every time step for a base and initial refinement level of two and a maximum refinement level of six.

The results of an experiment running three separate simulations for each configuration and reporting the minimum measured time can be found in Table~\ref{tab:numres:perf:performance}. The results on the uniformly refined mesh show good speedup up to four MPI ranks, with a parallel efficiency above 80\%, dropping below 40\% at 32 MPI ranks. The share of runtime spent in each of the two solvers stays essentially constant across all rank counts, indicating that both components scale alike. The proportion of time spent in the gravitational solver is comparable to other available solvers for the dynamics of self-gravitational gases. 

In the benchmark results, it can be seen that a significant amount of computational time with enabled \(h\)-adaptivity is spent on adapting the mesh. This is because we request a mesh adaptation after each time step, which sets up all the data structures for the multigrid solver anew. When solely comparing the plain times spent in each of the solvers, it is noticeable that the \(h\)-adaptive run is performed faster than in the non-adaptive run. Also, \(h\)-adaptive simulation runs require less memory than those using uniformly refined meshes, as the domain is only as finely resolved as it needs to be with respect to the mesh refinement indicator. For a more practical setup, a choice has to be made for how often indicators regarding mesh adaptation are to be checked and the mesh adapted afterwards. However, one should be aware that mesh refinement can also stabilize a numerical simulation and the interval in which the indicators are computed should not be chosen too large.

\begin{table}
    \centering

    uniformly refined meshes\\
    \begin{tabular}{ccccc}
        \toprule
         MPI ranks & \multicolumn{2}{c}{total time/speedup} & \multicolumn{1}{c}{hydrodynamics} & \multicolumn{1}{c}{gravity} \\
         \midrule
         \nph1 & 23.5\nph h &  --  & 33.7\% & 65.3\% \\
         \nph2 & 12.0\nph h & \nph1.96 & 34.0\% & 65.2\% \\
         \nph4 & \nph7.06 h & \nph3.33 & 34.7\% & 64.2\% \\
         \nph8 & \nph5.94 h & \nph3.96 & 32.3\% & 66.6\% \\
         16 & \nph3.98 h & \nph5.90 & 31.1\% & 67.8\% \\
         32 & \nph1.95 h & 12.1\nph & 32.3\% & 66.0\% \\
         \bottomrule
    \end{tabular}

    \vspace{\baselineskip}
    with adaptive mesh refinement\\
    \begin{tabular}{cccccc}
        \toprule
         MPI ranks & \multicolumn{2}{c}{total time/speedup} & \multicolumn{1}{c}{hydrodynamics} & \multicolumn{1}{c}{gravity} & \multicolumn{1}{c}{mesh adaptation} \\
         \midrule
         1 & 7.50 h &  --  & 16.5\% & 37.3\% & 45.6\% \\
         2 & 4.04 h & 1.86 & 16.0\% & 36.6\% & 46.6\% \\
         4 & 2.13 h & 3.52 & 15.9\% & 35.9\% & 47.3\% \\
         8 & 1.21 h & 6.20 & 16.2\% & 34.6\% & 48.0\% \\
         \bottomrule
    \end{tabular}
    \caption{Total time and speedup for different numbers of MPI ranks, along with the percentage of runtime spent in each component of the coupled Trixi.jl--deal.II solver, for the three-dimensional spherical Sedov blast wave.}
    \Description{In these tables we show the runtime and percentages of runtime spent in the subroutines of the hydrodynamical solver, the gravitational solver and runtime spent performing adaptive mesh refinement. We took benchmarks of a run on a uniformly refined mesh, for which the parallel efficiency decreases from 98\% on two MPI ranks to below 40\% on 32 MPI ranks. We also took benchmarks of a run with adaptive mesh refinement occurring every time step, which retains a parallel efficiency above 75\% up to eight MPI ranks. In both cases the number of MPI ranks used were taken to be powers of two.}
    \label{tab:numres:perf:performance}
\end{table}

\section{Conclusions \& Outlook} \label{sc:outlook}

In this work, we have proposed a modular and easily extensible framework for coupled numerical solvers of multi-physics problems based on a minimally invasive cross-language coupling of the Trixi.jl framework and the finite element library deal.II. As an application, a parallel adaptive coupled solver for the numerical simulation of the dynamics of Newtonian self-gravitational gases has been presented. We have shown convergence results and practical applicability on common test cases as well as strong scaling for MPI parallel executions, with a parallel efficiency above 75\% up to eight ranks in the adaptive simulation runs. The source code and scripts for reproducing the results presented in this work have been made publicly available at \cite{ehlert2026multisolverRepro}.

The current coupling keeps two copies of the mesh in the two solvers, each with full data structures, for which we faced challenges in keeping the representation synchronized. Therefore, future work will consider sharing a mesh and then updating a thin view of this mesh along with buffers and solver data structures due to changes in the connectivity on the deal.II side of the coupling, enabled by the \texttt{parallel::fullydistributed::Triangulation} class of deal.II. This alternative approach could lower both the setup cost and memory requirements. Furthermore, the complete re-computation of data structures on multigrid levels adds overhead, which we plan to avoid by identifying which parts of the previous data structures can be re-used.

While the present work offers a proof-of-concept of coupling between the two solvers, a multitude of features available in both frameworks have not yet been coupled. This includes more flexible boundary conditions, curved meshes, non-Cartesian meshes via T8code or loading meshes from files. Also, both Trixi.jl and deal.II support off-loading computations to GPUs. As a follow-up to this manuscript, we plan to publish a somewhat extended version of this coupling framework as an extension of Trixi.jl.


\section*{CRediT authorship contribution statement}
\textbf{Vivienne Ehlert:} Investigation, Software, Visualization, Writing -- original draft, Writing -- review \& editing;
\textbf{Gregor J.~Gassner:} Conceptualization, Funding Acquisition, Writing -- review \& editing;
\textbf{Martin Kronbichler:} Conceptualization, Software, Writing -- original draft, Writing -- review \& editing;
\textbf{Hendrik Ranocha:} Conceptualization, Software, Writing -- review \& editing;
\textbf{Michael Schlottke-Lakemper:} Conceptualization, Funding Acquisition, Software, Supervision, Writing -- review \& editing

\begin{acks}
Gregor J.~Gassner and Michael Schlottke-Lakemper acknowledge support by the DFG research
unit FOR 5409 “SNuBIC” (Structure-Preserving Numerical Methods for Bulk- and Interface
Coupling of Heterogeneous Models).
\end{acks}

\bibliographystyle{ACM-Reference-Format}
\bibliography{bibliography}

\end{document}